\documentclass[11pt]{article}	
\begin{document}
\center{\bf Singular Value Decomposition at FIFA 2022 \\
\it Peter Zizler and Mandana Sobhanzadeh\\
Mount Royal University,
Calgary}

\abstract{ The singular value decomposition is arguably one of the most fundamental results in linear algebra. While rigorous proof of this result is of importance, equally important
is the motivation in the applied settings. We provide a lively and quite intuitive presentation on the appearance of the singular value decomposition of a matrix 
in the round robin tournaments, as well as the polar decomposition, all in the context of FIFA 2022 World Cup. This exposition is intended to be implemented in a class setting or given as a take home project in the second linear algebra course.
Given the popularity of this world event we expect students to be interested and engaged throughout the analysis and ensuing conversations.
Discussions among students should be enriching as they attempt to analyse the various FIFA groups and give comparisons. Furthermore, ideas in regards how to best interpret the singular vectors should be quite interesting.
Basic notions and results from the first year linear algebra course are assumed.}
		
\section{Introduction} 
Our mathematical essay can be seen as a motivation to learn the need for the 
singular value decomposition of a matrix as well as its closely related result, the polar decomposition of a matrix. We will see how the singular value decompsotion of a matrix has a practical meaning in round robin tournaments. Moreover, students will also appreciate and discuss the 
closely related polar decomposition of a matrix, a result that will find a specific interpretation in the round robin tournament matrix as well. 

Given the results of the round robin tournament it is natural to seek the corresponding teams' offense and defense scores. As a result we have a motivation for the rank one matrix approximation, the tensor product
approximation, of the tournament matrix $A$. It is intuitive that the team's offense and defense ability can be estimated by goal counts, for and against, in round robin tournaments. Students will be encouraged to explore other options, possibly weighted goal counts, both for and agaisnt. The choice
of the appropriate wieights might stimulate enriching discusions that will enhance future deeper understadning the the forthcoming singular value decomposition of $A$, in particular the interpretation of its dominant singular vectors. 
Students will find the offense and defense scores for the teams, stemming out of the singular value decomposition of a matrix, are sensitive to goal scoring on good or bad defensive teams, and goals against from good or bad offense teams. 
In matrix theory language, we motivate the best least squares approximation of $A$ by a rank one matrix (tensor),

$$
A \approx \sigma_1 {\bf u}_1 {\bf v}_1^T
$$

\noindent where ${\bf  u}_1$ are the normalized team offense scores and ${\bf v}_1$ are the normalized team defense scores, $\sigma_1$ being a scaling factor. 

We continue to analyze all groups at FIFA 2022 using these techniques with further finer analysis of team interactions in the tournaments based on the other rank one approximations in the singular value decomposition. Connecting our exposition to a popular event like the FIFA 2022 World Cup can motivate the students to engage and sustain interest. On the other hand, these methods can be applied in any round robin tournament for any relevant sport activities should the students seek their favourite sports settings.

During the FIFA World Cup in Qatar, Fall $2022$, a very interesting scenario arose when Argentina, Poland, Mexico and Saudi Arabia played round robin games, from
which two teams would advance. This was the group C at the soccer tournament. If a team wins against another team the winning team gets awarded three points, the team that lost gets zero points. If the match results in a tie 
then both teams get awarded one point.
The two teams with the most points advance. However, many teams might have the same number of points, in fact, this happens relatively often.
In such an event, the goal difference then decides. One collects all the goals scored by the team (in the three matches) minus the goals conceded by the team (in the three matches).
If this still fails to decide, then only the number of goals scored by the team in the three matches decide without any consideration to the conceded goals by the team.
If this does not decide, then the team head to head results decide who advances. If this still does not decide then the teams with fewer disciplinary 
actions against them advance, a rule that involves yellow and red cards. This is referred to as the fair play rule. Finally, if this still does not decide, then a toss of a coin will.

Argentina had $6$ points, Poland had $4$ points, Mexico had $4$ points and Saudi Arabia had $3$ points. Argentina would advance, Saudi Arabia was in the last place. 
The decision had to be made between Poland and Mexico. Poland goal differential was $0$ where as Mexico goal diffential was $-1$. Poland advanced. 
However, the excitement at the time was the following. Argentia played Poland 
and Mexico played Saudi Arabia as the last two games in the group simultaneously. Argentina lead $2:0$ and Mexico lead $2:0$ to the very last minutes of the games. If the score
remained as such then the decision would had to be made by the fair play rule. However, Saudi Arabia scored.

The results of the round robin tournament reflect the teams' offensive ability as well as their defensive ability. As the teams play in the tournament they score goals against other teams during the matches. 
The total amount of goals scored by a team in the tournament is a possible way to measure the team's offensive ability, however it values each goal scored equally, goal against a good defensive team in the same way
as a goal against a poor defensive team. It is harder to score against defensive teams. Similarly, each team has a certain amount of goals scored on them. The count of goals against is a possible measure for the team's defensive ability, 
however, it is easier to defend against a poor offensive team as opposed to a good offensive team. Offense and defense scores for the teams, which are sensitive to this aspect, can be obtained from the singular
value decomposition of a matrix. The matrix in question is the offense, defense performance matrix associated with the given round robin tournament.   

We will assign normalized offense scores to each team in the tournament as well as the normalized defense scores for each team. This will in turn induce the offense and defense rankings for the teams in the tournament. Compared to the offense and defense rankings stemming out of total goals scored in the tournament, and goals scored against in the tournament, the singular value decomposition rankings break ties that often arise in the tournament. 
Many times the scores from the singular value decompostion rankings, offense and defense, yield the same ranking results as the total goals count, however there are cases, some groups at FIFA 2022, where the rankings obtained are different. We will show these results in our paper. We will analyze each group at the 2022 FIFA world cup, groups A through H. We first illustrate the technique on the group C. The techniques are applicable to any round robin tournament. The analysis drawn can be utilized to predict the team scoring outcomes should the teams play again under the same attributes. The ideas appearing in the paper are inspired by the work of \cite{James:2013}.
For more works on related topics see \cite{Hepler:2016}. We now present the singular value decomposition, the offense, defense decomposition for the {\bf group C} at FIFA 2022.

\section{\bf The Singular Value Decomposition} The round robin results for the {\bf group C} at FIFA 2022 were

\begin{center}
\begin{tabular}{ c c c c c}
  & Argentina & Poland & Mexico & Saudi Arabia \\ 
 Argentina & x & 2:0 & 2:0 & 1:2 \\  
 Poland & 0:2 & x & 0:0 & 2:0  \\
 Mexico & 0:2 & 0:0 & x & 2:1 \\ 
Saudi Arabia & 2:1 & 0:2 & 1:2 & x  \\ 
\end{tabular}
\end{center}

We consider the results in this tournament using the lenses of the team's offense and defense ability. What counts now are the goals 
scored by the teams and against, winning, tying or losing a game is irrelevant. To this end, we form the offense defense performance matrix $A$ reflecting the tournament results

$$
A=
\left( \begin{array}{cccc}
 1.1667  & 2.00 & 2.00  & 1.00\\
 0  & 0.6667  & 0 &  2.00\\
 0  & 0 &  0.8333 &  2.00\\
 2.00 & 0 &  1.00 &  1.3333\\
\end{array} \right)
$$

\noindent The off diagonal value $A(i,j)$, $i \neq j$, is the number of goals scored by the team $i$ against the team $j$. However, lively discussions in regards to the diagonal entires in the matrix $A$ should ensue among students. It reflects a hypothetical scenario where a team plays against itself and it should capture 
both the offense and the defensive ability of the team. One option is to set the diagonal value $A(i,i)$ to be the mean of all row $i$ and column $i$ entries, excluding the diagonal, all $6$ of them. 
The diagonal entry $A(i,i)$ thus averages out the goals scored by the team $i$ and goals scored against the team $i$. Having said that, we presume that students might suggest further ideas on this. We consider the singular value decomposition of $A$ 

\begin{eqnarray*}
A & = & UDV^T \\
& = & 
\left( \begin{array}{cccc}
  0.6652 &  0.5771 & -0.4710 &  0.0512 \\
  0.3629 & -0.5899 &  -0.2826 & -0.6636 \\
  0.4038 & -0.5547 & -0.0304 &  0.7269 \\
  0.5127 &  0.1057 &  0.8351 & -0.1692 \\
\end{array} \right) \\
& & 
\left( \begin{array}{cccc}
  4.3328    &    0    &    0   &     0 \\
        0  &  2.1135   &     0    &    0 \\
        0   &     0 &  1.5973   &     0 \\
        0    &    0  &      0  & 0.6971
\end{array} \right)
\left( \begin{array}{cccc}
0.4158  & 0.4186  & 0.7016 &  -0.3997 \\
  0.3629  & 0.3600 & -0.7078 & -0.4877 \\
  0.5030  & 0.3774  & -0.0828  & 0.7731 \\
  0.6652  & -0.7435  & 0.0102  & -0.0687 \\
\end{array} \right)
\end{eqnarray*}

\begin{eqnarray*}
A & = &  \sigma_1{\bf u}_1{\bf v}^T_1 + \sigma_2{\bf u}_2{\bf v}^T_2 + \sigma_3{\bf u}_3{\bf v}^T_3+ \sigma_4{\bf u}_4{\bf v}^T_4 \\
& = & A_1+A_2+A_3+A_4 \\
& = & 
4.3328\left( \begin{array}{c}
0.6652 \\
0.3629 \\
0.4038 \\
0.5127
\end{array} \right)
\left( \begin{array}{cccc}
0.4158 & 0.3629 & 0.5030 & 0.6652 \\
\end{array} \right) \\
& + & 
2.1135\left( \begin{array}{c}
  0.5771 \\
 -0.5899  \\
-0.5547   \\
0.1057 \\
\end{array} \right)
\left( \begin{array}{cccc}
0.4186  &  0.3600 &   0.3774 &  -0.7435 \\
\end{array} \right) \\
& + &
 1.5973\left( \begin{array}{c}
    -0.471010  \\
-0.282612  \\
-0.030440   \\
0.835077 \\
\end{array} \right)
\left( \begin{array}{cccc}
 0.701562  & -0.707705  & -0.082831 &   0.010201 \\
\end{array} \right) \\
& + &
 0.6971\left( \begin{array}{c}
  0.051229  \\
-0.663640   \\
0.726862  \\
-0.169204 \\
\end{array} \right)
\left( \begin{array}{cccc}
 -0.399693   & -0.487699   & 0.773091 &  -0.068748 \\
\end{array} \right) \\
\end{eqnarray*}

In particular,

\small
\begin{eqnarray*}
A & = & \left( \begin{array}{cccc}
1.1982 &  1.0458 &  1.4497 &  1.9170 \\
   0.6537 &  0.5705  & 0.7909  & 1.0458 \\
   0.7274 &  0.6348  & 0.8800  & 1.1637 \\
   0.9236 &  0.8061  & 1.1174  & 1.4776 \\
\end{array} \right)
+
\left( \begin{array}{cccc}
   0.5106 &  0.4392  & 0.4604  &-0.9069 \\
  -0.5219 & -0.4489&  -0.4706 &  0.9270 \\
  -0.4907 & -0.4221&  -0.4425 &  0.8716 \\
   0.0935 &  0.0804&   0.0843 & -0.1660 \\
\end{array} \right) \\
& + & 
\left( \begin{array}{cccc}
 -0.5278  & 0.5325 &   0.0623  &-0.0077 \\
  -0.3167&   0.3195  & 0.0374  &-0.0046 \\
  -0.0341  & 0.0344   &0.0040 & -0.0005 \\
   0.9358  &-0.9440  &-0.1105 &  0.0136 \\
\end{array} \right)
+
\left( \begin{array}{cccc}
 -0.0143 & -0.0174  & 0.0276 & -0.0025 \\
   0.1849 &  0.2256&  -0.3577  & 0.0318 \\
  -0.2025 & -0.2471&   0.3917 & -0.0348 \\
   0.0471  & 0.0575&  -0.0912 &  0.0081 \\
\end{array} \right) \\
\end{eqnarray*}
\normalsize

\noindent The vectors $\left\{ {\bf v}_ i\right\}^4_{i=1}$ form an orthonormal basis for ${\bf R}^4$, the same holds for the set  $\left\{ {\bf u}_ i\right\}^4_{i=1}$. We have $A{\bf v}_i = \sigma_i {\bf u}_i$ for $i \in \{1,2\ldots,4 \}$. The matrix $A_1$ is the best (in least squares sense)
rank one matrix approximating the matrix $A$. In particular, 

$$
{\rm min} \left\{ ||A-B||_2 \mbox{ where }  B \mbox{ is rank one matrix }\right\} = ||A-A_1||_2
$$

\noindent here $||C||^2_2$ denotes the square of the Frobenius norm of the matrix $C$, the sum of the squares of the absolute values of the entries in the matrix $C$. The matrix $A_2$ is the best (in least squares sense) rank one matrix approximating the matrix $A-A_1$. We continue this recursively and the procedure ends after four steps. For a very nice exposition of the singular value decomposition we refer the reader to \cite{Strang:2006}. The singular vectors ${\bf u}_1$ and ${\bf v}_1$ corresponing to the largest singular value $\sigma_1$ of $A$ are guaranteed to have nonnegative entries, due to 
the Perron - Frobenius Theorem, since the entries in the matrix $A^TA$ as well as $AA^T$ are all nonnegative. For a good reference on the Perron-Frobenius Theorem see \cite{Lancaster:1985}.

\subsection{\bf The Offense and Defense Scores} The (normalized) offense scores are the entries in the vector \\
${\bf u}_1 = (0.6652,   0.3629,   0.4038 ,  0.5127)^T$. In particular, Argentina is awarded the highest score of $0.6652$, then Saudi Arabia with score of $0.5127$, then 
Mexico with $0.4038$ and Poland last with $0.3629$. The higher the offense score the better the offense rating. The (normalized) 
defense scores are the entries in the vector ${\bf v}_1 = (0.4158,   0.3629 ,  0.5030 ,  0.6652)^T$, higher the score worse the defense.
In the defense scores Poland is the best with $0.3629$, then Argentina with $0.4158$, then Mexico with $0.5030$ and Saudi Arabia last with $0.6652$. We note that normalized here means the sum of the squares of the entries
in the vector is unity. The best rank one approximation for $A$ in the least squares sense is given by $A_1=\sigma_1{\bf u}_1{\bf v}^T_1$. The singular value decomposition scores yield the following ranking. Offense: Argentina, Saudi Arabia, Mexico, Poland. Defense: Poland, Argentina, Mexico, Saudi Arabia. We note the total goal count yields the same ranking, albeit induces ties. In particular Argentina scored $5$ goals, Poland $2$ goals, Mexico $2$ goals, and Saudi Arabia 
$3$ goals. We define the goal scored vector $(5,2,2,3)^T$ and upon normalization obtain the vector $(0.7715,   0.3086,   0.3086,   0.4629)^T$. This vector can be seen as a crude approximation of the singular vector ${\bf u}_1$.
Argentina allowed $2$ goals, Poland $2$ goals, Mexico $3$ goals, and Saudi Arabia 
$5$ goals. We define the goal allowed vector $(2,2,3,5)^T$ and upon normalization obtain the vector $(0.3086 ,  0.3086,   0.4629,   0.7715)^T$. This vector can be seen as a crude approximation of the singular vector ${\bf v}_1$.
For more work on finer approximations of the first singular vectors, that involves weighted row and column sums, we refer the reader to \cite{Zizler:2020}.

If the results of the tournament where only determined by the teams offense and defense scores, then the tournament performance matrix would be

$$
A_1 =\left( \begin{array}{cccc}
 1.1984 &  1.0459 &  1.4497 &  1.9172 \\
   0.6538 &  0.5706 &  0.7909 &  1.0459 \\
   0.7275 &  0.6349&   0.8800  & 1.1638 \\
   0.9237 &  0.8062 &  1.1174 &  1.4777 \\
\end{array} \right)
$$

\noindent For example, the score of $A_1(2,3)=0.7909$ means that Poland should score $0.7909$ goals against Mexico. Similarly,
the score of $A_1(4,1)=0.9237$ means that Saudi Arabia should score $0.9237$ goals against Argentina. The diagonal entries are of interest as well. For example, the diagonal entry $(2,2)$ equals to $0.5706$. It indicates, based on the prediction
from offense and defense scores, how Poland would fare in offense and defense simultaneously. Do note however, that the diagonal entry at $(i,i)$ in the matrix $A_1$
is no longer the mean of the off diagonal entries in the row $i$ and column $i$.

The amount of goal scoring during the tournament that is explained by the teams' offense and defense scores is

$$
\frac{\sigma_1^2} {\sigma_1^2+\sigma_2^2+\sigma_3^2+\sigma_4^2} = 0.7144
$$

\noindent indicating that around $71$ \% of goals scored during the tournament can be explained by the teams offense and defense scores. 

\subsection{\bf The Correction Scores} Giving interpretation to the second order singular vectors in the decompostion should be quite interesting. We expect students to give ample feedback here as the interpretations are not immediate. Consider the (goal) correction matrix

$$
A-A_1=
\left( \begin{array}{cccc}
 -0.0317 &  0.9541  &  0.5503  &  -0.9172 \\
  -0.6538 &  0.0961 &  -0.7909  &   0.9541 \\
-0.7275 & -0.6349 &  -0.0467 &   0.8362 \\
   1.0763  & -0.8062 &   -0.1174 &  -0.1444 \\
\end{array} \right)
$$

\noindent with the understanding of the entries as follows. The score of $0.9541$, Argentina vs Poland, indicates that when Argentina actaully played Poland, the number of goals that Argentina scored against Poland 
was $0.9541$ higher than the offense and defense scores would predict. Similarly, the score of $-0.1174$, Saudi Arabia vs Mexico, indicates that when Saudi Arabia actaully played Mexico, the number of goals that Saudi Arabia
scored against Mexico was $0.1174$ lower than the offense and defense scores would predict. The diagonal entry at $(4,4)$ has a value of $-0.1444$. It indicates that the actual amount of goals Saudi Arabia scored and allowed  
was lower than was predicted by the offense and defense scores. We now give interpretation to the vectors 

$$
{\bf u}_2=(0.5771, -0.5899,  -0.5547,   0.1057 )^T \mbox{ ; } {\bf v}_2=(0.4186,    0.3600,    0.3774,   -0.7435 )^T
$$
 
\noindent We are approximating, in the least squares sense, the correction matrix $A-A_1$ by a rank one matrix 

\begin{eqnarray*}
A_2& = & \sigma_2{\bf u}_2{\bf v}^T_2 \\
& =  &
\left( \begin{array}{cccc}
   0.5106 &  0.4392  & 0.4604  &-0.9069 \\
  -0.5219 & -0.4489&  -0.4706 &  0.9270 \\
  -0.4907 & -0.4221&  -0.4425 &  0.8716 \\
   0.0935 &  0.0804&   0.0843 & -0.1660 \\
\end{array} \right)
\end{eqnarray*}

The matrix $A_2$ can be seen as a (rank one) predicted correction matrix. It predicts the number of goals to be scored by a team during a match as a correction to the expected number of goals to be scored based on the offense defense scores of the teams. 
The matrix $A_2$ is the best least squares rank one correction matrix, noting that the actual goals scored during the game involve two higher level correction matrices $A_3$, and $A_4$.  
The entries in the vector ${\bf u}_2$ and ${\bf v}_2$ can be thought of as a correction inducing scores, offense and defense respectively. Positive scores reinforce, negative scores oppose. 
We group the positive offense correction scores, ranked by magnitude, \{Argentina, Saudi Arabia\} and group the negative offense correction scores, ranked by absolute value of the magnitude, \{Poland, Mexico\}. 

Similarly, we group the positive defense correction scores, ranked by magnitude, \{Argentina, Mexico, Poland\} and detect the negative defense correction score for \{Saudi Arabia\}. When a team $i$ in the positive offense correction group 
plays a team $j$ in the positive defense correction group we expect the team $i$ to score more goals than expected (based on the offense defense scores) against the team $j$, proportional to magnitude of the respective correction scores. When a team $i$ in the negative offense correction group plays a team $j$ in the negative defense correction group we also expect the team $i$ to score more goals than expected against the team $j$, proportional to absolute values of the respective correction scores. 
However, when a team $i$ in the positive offense correction group plays a team $j$ in the negative defense correction group we expect the team $i$ to score fewer goals than expected against the team $j$, proportional to absolute values of the respective correction scores. Identical situation is detected for the negative offense correction group and the positive defense correction group combination.

The diagonal entries are of interest as well. The diagonal entry $(2,2)$ equals to $-0.4489$, correponds to team Poland. It indicates the team Poland is expected to correct with fewer goals scored and allowed than the offense defense team scores
would indicate. Loosely said, Poland is predicted to have less goal action around them than predicted based on offense defense ability. On the other hand, the diagonal entry $(1,1)$ equals to $0.5106$, corresponds to team Argentina. 
Argentina is predicted to have more goal action around them than predicted based on offense defense ability. We can assert, on some level, that the teams Poland, Mexico and Saudi Arabia were predicted to be more boring than expected (based on the offense defense skills). Ranking of predicted boredom is given by the absolute value of the diagonal entry, making Poland the most boring predicted and Saudi Arabia the least. Argentina is predicted to be more exciting.

We reiterate that the matrix $A_2$ is just the best rank one prediction for the actual correction matrix $A-A_1$, based on the defense and offense correction scores for the teams. The actual results in the tournament
are different as the matrices $A_3$ and $A_4$ were not involved. For example, the predicted outcomes for the matches of Saudi Arabia using its offense correction score were substantially different than the correction matrix $A-A_1$ actually corrected. In particular, the match Saudi Arabia vs Mexico was predicted by their correction scores to correct with more goals for Saudi Arabia, when in fact the actual result in the tournament was the opposite. 

As mentioned above the expected results drawn from the offense and defense team scores and the correction offense and defense scores can be utilized for the scenario where the teams play again in a round robin tournament under the 
same attributes. The techniques discussed can be applied to any round robin tournament of any size in any relevant sport settings. 

\section{Polar Decomposition}
Polar decompostion of a matrix is a closely related result to the singular value decomposition of a matrix, albeit its applications are not seen as often. Having said that, the polar decomposition has a rightful place in the context of round robin tournament matrices. We believe students 
of linear algebra will find this place for the polar decomposition quite intuitive, expecting many what if scenarios.
Let us return to the results of the {\bf group C}. The teams involved had their offense attributes captured by the singular vectors $\left\{{\bf u}_1,{\bf u}_2,{\bf u}_3,{\bf u}_4  \right\}$
and their defense attributes captured by $\left\{{\bf v}_1,{\bf v}_2,{\bf v}_3,{\bf v}_4  \right\}$. We now address the hypothetical tournament outcome where the teams' defense attributes 
are identical to those of their offense attributes. In particular, ${\bf v}_i={\bf u}_i$ for all $i \in \left\{1,2,3,4 \right\}$. The condition ${\bf v}_1={\bf u}_1$ would force the defense ability scores for the teams to match their offense ability scores. This would mean that a team that has a good offense ability would in turn 
suffer on defense ability and vice versa. High score on offense translates to the same high score on defense which indicates a poor defense. We expect lively discussions among students in regards to this, how to interpret the forthcoming symmetric tournament matrices in the polar decomposition.
The singular value decomposition of the matrix $A$ can be written in a polar decomposition of the matrix $A$ 

$$
A = U D V^T = \left(UDU^T\right)\left( UV^T \right) 
$$  

We set $P=UDU^T$ and the entries in the matrix $P$ provide the answer to our hypothetical tournament scenario. The matrix $P$ is a symmetric matrix, positive definite. We have

$$
P=
\left( \begin{array}{cccc}
    2.9772  &  0.5151 &   0.5359  &  0.9721 \\
    0.5151  &  1.7407  &  1.0040  &  0.3757  \\
    0.5359  &  1.0040  &  1.7266  &  0.6467 \\
    0.9721  &  0.3757  &  0.6467  &  2.2964 \\
\end{array} \right)
$$

\noindent Observe $||P||_2=||A||_2$. The matrix $P$ gives the performance results for the tournament under the hypothetical setting where the 
teams' defense attributes reflect those of their offense attributes. Each match between a team $i$ and team $j$ would end up in a draw. For example, if Poland played Saudi Arabia we would expect a score of $0.3757:0.3757$. Similarly, if Argentina played Mexico we would
expect a score of $0.5359:0.5359$. 

Equivalently, we now address the hypothetical tournament outcome where the teams' offense attributes 
reflect those of their defense attributes. Now we have ${\bf u}_i={\bf v}_i$ for all $i \in \left\{1,2,3,4 \right\}$. In particular, ${\bf u}_1={\bf v}_1$ would force the offense ability scores for the teams to match their defense ability scores. This would mean that a team that has a bad defense would have good offense and vice versa. We have the following

$$
A = U D V^T = \left( UV^T \right) \left(VDV^T\right)
$$  

We set $Q=VDV^T$ and the entries in the matrix $Q$ provide the answer to our hypothetical tournament scenario.  The matrix $Q$ is a symmetric matrix, positive definite. We observe

$$
Q=
\left( \begin{array}{cccc}
    2.0168  &  0.3150  &  0.9318  &  0.5711 \\
    0.3150  &  1.8103  &  0.9089  &  0.4919 \\
    0.9318  &  0.9089  &  1.8250  &  0.8183 \\
    0.5711  &  0.4919  &  0.8183  &  3.0887 \\
\end{array} \right)
$$

\noindent We also have $||Q||_2=||A||_2$. The matrix $Q$ gives the performance results for the tournament under this hypothetical setting where the 
teams offense attributes reflect those of their defense attributes. Once again, each match between a team $i$ and team $j$ would end up in a draw. For example, if Poland played Saudi Arabia in this setting we would expect a score of $0.4919:0.4919$. Similarly, if Argentina played Mexico in this setting we would
expect a score of $0.9318:0.9318$. 

Note that the diagonal entries in both $P$ and $Q$ are no longer the means of the off diagonal entries in the corresponding rows and columns. For more on a polar decomposition of a matrix see \cite{Lancaster:1985}.   

\section{FIFA Groups}
The draw procedure used in the $2022$ FIFA World Cup draw on the $1$st of April $2022$ is described here. The $32$ teams were divided into four pots based on the FIFA World Ranking
announced on the $31$st of March $2022$. The pot $1$ contained the host Qatar (assigned to group A) and the seven highest-ranked teams. The pot $2$ contained the teams ranked $8$th to $15$th and the
pot $3$ included the teams ranked $16$th to $23$rd, while the pot $4$ contained the teams ranked $24$th to $28$th and the two placeholders from the two inter-confederation play-offs and the winner of the UEFA play-off Path A. The draw sequence started with the pot $1$ and ended with pot $4$. Each pot was emptied before moving on to the next pot. Some draw conditions ensured geographic separation. For instance, no group can have more than one team from any continental confederation except for UEFA (AFC, CAF, CONMEBOL, CONCACAF) and each group should consist of at least one but no more than two European teams, see \cite{Csato:2022}.

As a result the World Cup had $8$ groups of $4$ teams. The top $2$ teams from each group advanced to the knockout stage. We now perform the singular value decomposition for all the round robin results for all the groups.

{\bf Group A}

\begin{center}
\begin{tabular}{ c c c c c}
  & Netherlands & Senegal & Ecuador & Qatar \\ 
 Netherlands & x & 2:0 & 1:1 & 2:0 \\  
 Senegal & 0:2 & x & 2:1 & 3:1  \\
 Ecuador & 1:1 & 1:2 & x & 2:0 \\ 
Qatar & 0:2 & 1:3 & 0:2 & x  \\ 
\end{tabular}
\end{center}

$$
A=
\left( \begin{array}{cccc}
1   & 2.00 & 1.00  & 2.00\\
 0.00  & 1.5  & 2.00 &  3.00\\
 1.00  & 1.00 & 1.6667  &  2.00\\
 0.00 & 1.00 &  0.00 & 1.3333  \\
\end{array} \right)
$$

\noindent Explained by offense and defense scores:  $0.9302$

$$
{\bf u}_1=( 0.5200,   0.6534 ,  0.4909 ,  0.2486)^T \mbox{ ; } {\bf v}_1=(0.1728  , 0.4716  , 0.4520 ,  0.7372)^T
$$

$$
{\bf u}_2=( 0.6962,  -0.5597 , -0.1970 ,  0.4039)^T \mbox{ ; } {\bf v}_2=( 0.4202 ,  0.6395,  -0.6326  ,-0.1197)^T
$$

\noindent The singular value decomposition scores yield the following ranking. Offense: Senegal, Netherlands, Ecuador, Qatar. Defense: Netherlands, Ecuador,Senegal, Qatar. 
The total goal count yields the same ranking, albeit induces ties.

{\bf Group B}

\begin{center}
\begin{tabular}{ c c c c c}
  & England & USA & Iran & Wales \\ 
 England & x & 0:0 & 6:2 & 3:0 \\  
 USA & 0:0 & x & 1:0 & 1:1  \\
 Iran & 2:6 & 0:1 & x & 2:0 \\ 
Wales & 0:3 & 1:1 & 0:2 & x  \\ 
\end{tabular}
\end{center}

$$
A=
\left( \begin{array}{cccc}
 1.8333  & 0.00 & 6.00  & 3.00\\
 0.00  & 0.5  & 1.00 &  1.00\\
 2.00  & 0.00 & 1.8333  &  2.00\\
 0.00 & 1.00 &  0.00 & 1.1667  \\
\end{array} \right)
$$

\noindent Explained by offense and defense scores: $0.9266$

$$
{\bf u}_1=(  0.8968 ,  0.1698  , 0.4013 ,  0.0761)^T \mbox{ ; } {\bf v}_1=(0.3169 ,  0.0209 ,  0.8142 ,  0.4859)^T
$$

$$
{\bf u}_2=( -0.3400 ,  0.2428 ,  0.5151 ,  0.7485 )^T \mbox{ ; } {\bf v}_2=(  0.2382,   0.5092 , -0.4992 ,  0.6593 )^T
$$

\noindent The singular value decomposition scores yield the following ranking. Offense: England, Iran, USA, Wales . Defense: USA, England, Wales, Iran. 
The total goal count yields the same ranking.

{\bf Group D}

\begin{center}
\begin{tabular}{ c c c c c}
  & France & Australia & Tunisia & Denmark \\ 
 France & x & 4:1 & 0:1 & 2:1 \\  
 Australia & 1:4 & x & 1:0 & 1:0  \\
 Tunisia & 1:0 & 0:1 & x & 0:0 \\ 
Denmark & 1:2 & 0:1 & 0:0 & x  \\ 
\end{tabular}
\end{center}

$$
A=
\left( \begin{array}{cccc}
1.5   & 4.00 & 0.00  & 2.00\\
 1.00  & 1.1667  & 1.00 &  1.00\\
 1.00  & 0.00 & 0.3333  &  0.00\\
 1.00 & 0.00 &  0.00 & 0.6667   \\
\end{array} \right)
$$

\noindent Explained by offense and defense scores: $0.8918$

$$
{\bf u}_1=(  0.9191 ,  0.3619 ,  0.0801 ,  0.1338 )^T \mbox{ ; } {\bf v}_1=( 0.3832 ,  0.8036 ,  0.0762 ,  0.4489)^T
$$

$$
{\bf u}_2=( 0.3303 , -0.5100,  -0.5812 , -0.5412 )^T \mbox{ ; } {\bf v}_2=( -0.7402,   0.4728,  -0.4581,  -0.1368)^T
$$

\noindent The singular value decomposition scores yield the following ranking. Offense: France,Australia, Denmark, Tunisia. Defense: Tunisia, France, Denmark, Australia. 
The total goal count yields the same ranking, albeit induces ties.

{\bf Group E}

\begin{center}
\begin{tabular}{ c c c c c}
  & Japan & Spain & Germany & Costa Rica \\ 
 Japan & x & 2:1 & 2:1 & 0:1 \\  
 Spain & 1:2 & x & 1:1 & 7:0  \\
 Germany & 1:2 & 1:1 & x & 4:2 \\ 
Costa Rica & 1:0 & 0:7 & 2:4 & x  \\ 
\end{tabular}
\end{center}

$$
A=
\left( \begin{array}{cccc}
 1.1667  & 2.00 & 2.00  & 0.00\\
 1.00  & 2  & 1.00 &  7.00\\
 1.00  & 1.00 & 1.8333  &  4.00\\
 1.00 & 0.00 &  2.00 &  2.3333  \\
\end{array} \right)
$$

\noindent Explained by offense and defense scores: $0.8782$

$$
{\bf u}_1=( 0.1428  , 0.7958  , 0.4989 ,  0.3122  )^T \mbox{ ; } {\bf v}_1=(0.1930 ,  0.2586 ,  0.2851  , 0.9026 )^T
$$

$$
{\bf u}_2=( 0.8832 , -0.3527 ,  0.1359 ,  0.2777 )^T \mbox{ ; } {\bf v}_2=( 0.3603 ,  0.3952 , 0.7323 , -0.4216)^T
$$

\noindent The singular value decomposition scores yield the following ranking. Offense: Spain, Germany, Costa Rica, Japan. Defense: Japan,Spain, Germany, Costa Rica.
The total goal count yields different offense ranking. Offense: Spain, Germany, Japan, Costa Rica. Defense ranking is the same.

{\bf Group F}

\begin{center}
\begin{tabular}{ c c c c c}
  & Morocco & Croatia & Belgium & Canada \\ 
 Morocco & x & 0:0 & 2:0 & 2:1 \\  
 Croatia & 0:0 & x & 0:0 & 4:1  \\
 Belgium & 0:2 & 0:0 & x & 1:0 \\ 
Canada & 1:2 & 1:4 & 0:1 & x  \\ 
\end{tabular}
\end{center}

$$
A=
\left( \begin{array}{cccc}
0.8333   & 0.00 & 2.00  & 2.00\\
 0.00  & 0.8333  & 0.00 &  4.00\\
 0.00  & 0.00 & 0.5  &  1.00\\
 1.00 & 1.00 &  0.00 &  1.5  \\
\end{array} \right)
$$

\noindent Explained by offense and defense scores: $0.8334$

$$
{\bf u}_1=(0.4799,   0.7784,   0.2072,   0.3476 )^T \mbox{ ; } {\bf v}_1=(0.1473 ,  0.1964,   0.2096 ,  0.9465 )^T
$$

$$
{\bf u}_2=( 0.8467 , -0.4965 ,  0.1338 , -0.1369 )^T \mbox{ ; } {\bf v}_2=(0.2895,  -0.2803 , 0.8962,  -0.1854 )^T
$$

\noindent The singular value decomposition scores yield the following ranking. Offense: Croatia, Morocco, Canada, Belgium. Defense: Morocco, Croatia, Belgium, Canada.
The total goal count yields the same ranking, albeit induces ties.

{\bf Group G}

\begin{center}
\begin{tabular}{ c c c c c}
  & Brazil & Switzerland & Cameroon & Serbia \\ 
 Brazil & x & 1:0 & 0:1 & 2:0 \\  
 Switzerland & 0:1 & x & 1:0 & 3:2  \\
 Cameroon & 1:0 & 0:1 & x & 3:3 \\ 
Serbia & 0:2 & 2:3 & 3:3 & x  \\ 
\end{tabular}
\end{center}

$$
A=
\left( \begin{array}{cccc}
0.6667   & 1.00 & 0.00  & 2.00\\
 0.00  &  1.1667 & 1.00 &  3.00\\
 1.00  & 0.00 &  1.3333 &  3.00\\
 0.00 & 2.00 &  3.00 & 2.1667   \\
\end{array} \right)
$$

\noindent Explained by offense and defense scores: $0.8565$

$$
{\bf u}_1=(0.3203 ,  0.5207,  0.5007 ,  0.6129 )^T \mbox{ ; } {\bf v}_1=(0.1135,   0.3420 ,  0.4808 , 0.7994 )^T
$$

$$
{\bf u}_2=(0.3974  , 0.2011 ,  0.4696,  -0.7623  )^T \mbox{ ; } {\bf v}_2=(0.3352,  -0.4073,  -0.6661 ,  0.5273 )^T
$$

\noindent The singular value decomposition scores yield the following ranking. Offense: Serbia, Switzerland, Cameroon, Brazil. Defense: Brazil, Switzerland, Cameroon, Serbia .
The total goal count yields the same ranking, albeit induces ties.

{\bf Group H}

\begin{center}
\begin{tabular}{ c c c c c}
  & Portugal & South Korea & Uruguay & Ghana \\ 
 Portugal & x & 1:2 & 2:0 & 3:2 \\  
 South Korea & 2:1 & x & 0:0 & 2:3  \\
 Uruguay & 0:2 & 0:0 & x & 2:0 \\ 
Ghana & 2:3 & 3:2 & 0:2 & x  \\ 
\end{tabular}
\end{center}

$$
A=
\left( \begin{array}{cccc}
 1.1667  & 1.00 & 2.00  & 3.00\\
 2.00  & 1.3333  & 0.00 &  2.00\\
 0.00  & 0.00 &  0.6667 &  2.00\\
 2.00 & 3.00 &  0.00 &  2.0000\\
\end{array} \right)
$$

\noindent Explained by offense and defense scores: $0.8309$

$$
{\bf u}_1 =(0.5752,    0.4809 ,   0.2508 ,   0.6124)^T \mbox{ ; } {\bf v}_1=(0.4593 ,   0.4908  ,  0.2118   , 0.7094)^T
$$

$$
{\bf u}_2 =(0.5954,   -0.2229 ,   0.4993 ,  -0.5887)^T \mbox{ ; } {\bf v}_2 =(-0.3590,   -0.5676  ,  0.5892  ,  0.4492)^T
$$

\noindent The singular value decomposition scores yield the following ranking. Offense: Ghana, Portugal, South Korea, Uruguay . Defense: Uruguay, Portugal, South Korea, Ghana.
The total goal count yields different offense ranking. Offense: Portugal. Ghana, South Korea, Uruguay. Defense ranking is the same.

\end{document}